\input amstex
\documentstyle{amsppt}
\document
\magnification=1200
\NoBlackBoxes
\vsize17cm
\hoffset=0.7in
\voffset=1in
\nologo
\def\P{\bold{P}}
\def\Q{\bold{Q}}
\def\R{\bold{R}}
\def\C{\bold{C}}
\def\Z{\bold{Z}}


\bigskip

\centerline{\bf REMARKS ON MODULAR SYMBOLS}

\medskip

\centerline{\bf  FOR MAASS WAVE FORMS}

\medskip

\centerline{\bf Yuri I. Manin}

\medskip

\centerline{ Max--Planck--Institut f\"ur Mathematik, Bonn, Germany,}

\centerline{and Northwestern University, Evanston, USA}

\bigskip

\hfill{\it To Professor F.~Hirzebruch, with friendship and admiration,}

\hfill{\it for his anniversary.}

\medskip

{\bf Abstract.}  In this paper I introduce modular symbols
for Maass wave cusp forms. They appear in the guise of
finitely additive functions on the Boolean algebra
generated by intervals with non--positive rational ends,
with values in analytic functions
(pseudo--measures in the sense of [MaMar2]).
After explaining the basic issues and analogies
in the extended Introduction, I construct modular symbols
in the sec. 1 and the related L\'evy--Mellin transforms in the sec. 2.
The whole paper is an extended footnote to the
Lewis--Zagier fundamental study [LZ2].

\bigskip

\centerline{\bf \S 0. Introduction}

\medskip

{\bf 0.1. Period polynomials and period functions.} Let
$u(\tau )\in S_{2k}(SL(2,\Z))$ be a cusp form of an integer {\it  weight} 
$2k>0$
for the full modular group. This means that it is holomorphic
in the upper half plane, the tensor $u(\tau ) (d\tau)^k$ is $SL(2,\Z))$--invariant, and
$u(\tau )$ vanishes at cusps.

\smallskip

Its {\it period polynomial} is defined as the following integral:
$$
\psi (z)=\psi_u(z):=\int_0^{i\infty} u(\tau )(z-\tau )^{2k-2}d\tau
\eqno(0.1)
$$
Here $z$ is, for the time being, an auxiliary formal variable.

\smallskip
 
 One remarkable discovery in the theory of modular
functions was a possibility to develop its versions for a certain set
of {\it complex weights} $2s$ (replacing the former $2-2k$). 
This spectrum consists of the (doubled) zeroes of Selberg's
zeta function $Z(s)$ of $SL(2,\Z )$ acting on the upper half--plane, 
or equivalently, those
values of $s$ for which the Mayer transfer operator $\Cal{L}^2_s$
(cf. [May1], [May2])
has 1 as its eigenvalue: see  [LZ1] for a short review and
[LZ2] for a comprehensive exposition.

\medskip

{\bf 0.2. Classical modular symbols.} The classical modular symbols
of weight $2k$ for $SL(2,\Z)$, in one of their guises, can be defined simply
as integrals 
$$
\int_{\alpha}^{\beta} u(\tau )(z-\tau )^{2k-2}d\tau
\eqno(0.2)
$$
where this time $\alpha ,\beta \in \P^1(\Q)$ are
arbitrary cusps, and the integration is taken along, say,
the hyperbolic geodesic connecting $\beta$ to $\alpha .$

\smallskip

More precisely, the modular symbol $\{\alpha ,\beta\}_k$
(for the full modular group) is the integral (0.2)
considered as a linear map 
$$
\{\alpha ,\beta\}_k :\ S_{2k}(SL(2,\Z)) \to \C[z] .
\eqno(0.3)
$$
In the next subsections, we will briefly recall the number--theoretic motivations
for considering (0.3). A geometric interpretation
of (0.3), after a dualization, runs as follows: this integral
expresses the pairing between the Hodge cohomology
and the Betti homology of the moduli space
$\overline{M}_{1,2k-2}$ of elliptic curves with marked points
(cf. [Sh1], [Sh2] for a version involving Kuga varieties rather
than moduli spaces). 

\smallskip

The modular symbols (0.3) satisfy the following
simple functional equations:
$$
\{\alpha ,\beta\}_k + \{\beta ,\gamma \}_k +\{\gamma ,\alpha \}_k=0,\
\{\alpha ,\beta\}_k + \{\beta ,\alpha\}_k= \{\alpha ,\alpha \}_k=0 ,
\eqno (0.4)
$$
Thus they can be extended to a $\C[z]$--valued finitely additive function on the Boolean
algebra generated by (positively oriented) segments with rational ends
in $\P^1(\R)$. We sometimes call such a function {\it a pseudo--measure}, as in [MaMar2].
The variable change formula applied to (0.2) leads to an additional property
of this particular pseudo--measure, which we call
its {\it modularity}:
$$
\{g(\alpha ), g(\beta )\}_k =g \{\alpha ,\beta\}_k .
\eqno(0.5)
$$
Here $g\in  SL(2,\Z)$  acts on $\P^1(\Q)$ by fractional linear
transformations, and on polynomials of degree $\le 2k-2$
by a natural twisted action.

\medskip

A pseudo--measure can in principle take values in any abelian group,
and modularity condition (0.5) makes sense if this group is a left $SL(2,\Z )$--module.
If the group of values has no 2-- and 3--torsion, the last two equations in
(0.4) follow from the first one.

\medskip

{\bf 0.3. Modular symbols for Maass cusp forms.} The first goal of this 
note is to extend the definition of
$\{\alpha ,\beta\}_k$ to complex weights for which there
exist non--trivial Maass cusp forms. We take the formula (0.2)
as our starting point and look for its analogs in the Lewis--Zagier theory.
We are interested mostly in complex critical zeroes/weights
for which $\roman{Re}\,s=\dfrac{1}{2}.$

\smallskip

Tracing parallels with the classical theory, one should keep
in mind that certain classical objects have more (or less)
than one parallel in the new setting. 

\smallskip

For example,
the most straightforward analogs of  $u(\tau )\in S_{2k}(SL(2,\Z ))$
apparently are the Maass wave cusp forms, introduced in [M], -- 
smooth $SL(2,\Z)$--invariant
functions on $H$ satisfying the hyperbolic
Laplace equation $\Delta u= s(1-s) u$ and certain
growth/vanishing conditions. An appropriate version 
of the period polynomial (0.1) for such a form is its period function
$\psi_u(z)$, this time a holomorphic function
of our former auxiliary variable $z$.

\smallskip

However, the relationship between $u$ and $\psi_u$
as it is first explained in sec. 1 of Ch. I of [LZ2], does not look
at all like (0.1) and passes through  three
 intermediate steps 
$u \leftrightarrow L_{\varepsilon} \leftrightarrow f \leftrightarrow  \psi$.

\smallskip

To the contrary, the structure of (0.1) is reproduced in the formula
 $$
\psi (z)=\int_{-\infty}^0 (z-t)^{-2s} U(t)dt .
\eqno(0.6) 
$$
(see [LZ2], page 221), in which $U(t)dt$ denotes a certain distribution
on $\R$, called ``the boundary value'' of $u(\tau )$. Therefore,
it is this distribution that in our context seems to be a more adequate
analog of a classical cusp form, the more so
that its $SL(2,\Z)$--invariance property involves
an explicitly weighted action of the modular group,
$$
U\left(\frac{at+b}{ct+d}\right)= |ct+d|^{2-2s}  U(t).
\eqno(0.7)
$$
whereas a Maass form is simply $SL(2,\Z )$--invariant.

\smallskip

 The formula (0.6) seems to offer
a straightforward way to do it -- just consider the integrals
$$
\int_{\alpha}^{\beta} (z-t)^{-2s} U(t)dt .
$$
Formal manipulations with such integrals  are simple and
seemingly prove (0.4) and (0.5), and we reproduce them
for their heuristic value.
However, these calculations cannot be taken literally, because
the characteristic functions of the intervals with rational ends
do not belong to the space of test functions for the 
distribution $U(t)$. 

\smallskip

Thus we have to find a way around this difficulty.

\smallskip

In fact, there are at least two different ways. One of them starts with
the three--term functional equation for the period function $\psi (z)$,
proceeds with pure algebra,
and works also for Lewis--Zagier's ``period--like functions''.

\smallskip

Another method is applicable only to
the period functions of Maass forms $u$ and 
uses the Lewis--Zagier formula of the form
$$
\psi (z ) =\int^0_{-\infty} \{u, R^s_z\} (\tau )
$$
where the integrand is a closed 1--form depending on
$z$ as a parameter (its structure is described in the main text below). 
One can then integrate this form along a path that may this
time connect two arbitrary
cusps, thus getting another analog of (0.2).

\smallskip

These two constructions form the content of sec. 1 below.

\medskip

{\bf 0.4.  Mellin transform and classical  modular symbols.} 
Now we will explain some of our motivations.

\smallskip

Briefly, we want to describe a construction presenting
the Maass Dirichlet series as an integral
over, say, $[0,1/2]$, formally similar close to the
Mazur--Mellin transform in the theory of 
$p$--adic interpolation. We call such a representation
{\it the $\infty$--adic L\'evy--Mellin transform}, cf [MaMar2].
The integration measure in both cases is constructed
out of modular symbols.

\smallskip

Here is a sketch of the classical $p$--adic constructions.

\smallskip

The classical
theory of modular symbols as it was presented in [Ma1], [Ma2],
started with the following observations. Suppose that we are interested
in the calculation of some values (say, at integer points $\rho$)
of a Dirichlet series 
$$
L_{\kappa}(\rho )= \sum_{n=1}^{\infty} a_n\kappa (n)n^{-\rho}.
\eqno(0.8)
$$
where $(a_n)$ is a certain ``arithmetic'' function, and $\kappa$ is an additive
character of $\Z$ of finite order.
In the standard approach one first  introduces the Fourier series
$$
u_{\kappa}(\tau ):=\sum_{n=1}^{\infty} a_n\kappa (n)e^{2\pi in\tau}
\eqno(0.9)
$$
and then works with the Mellin transform
$$
\Lambda_{\kappa} (\rho ):=\int_0^{i\infty} u_{\kappa}(\tau)\left(\frac{\tau}{i}\right)^{\rho -1}
d\tau.
\eqno(0.10)
$$
which is related to (0.8) by the simple formula
$\Lambda (\rho)=i(2\pi )^{-\rho}\Gamma (\rho )L(\rho ).$ 

\medskip

Now, let $u(\tau ):=u_{\kappa_0}(\tau)$ where $\kappa_0$ is
identically $1$. Clearly, $u_\kappa (\tau )=u(\tau +\alpha )$
for a rational number $\alpha$ such that $\kappa (n) = e^{2\pi i\alpha n}$, so that we can write,
shifting the integration path,
$$
\Lambda_{\kappa}(\rho ):=\int_\alpha^{i\infty} u_{\kappa}(\tau)
\left(\frac{\tau -\alpha}{i}\right)^{\rho -1} d\tau
\eqno(0.11)
$$
Thus, if $\rho \ge 1$ is an integer, varying $\kappa$ in (0.8)
reduces to replacing $\tau^{\rho -1}$ in (0.10) by an arbitrary
polynomial of degree $\le \rho -1$ and allowing the
integration paths $(\alpha ,i\infty )$ with an arbitrary rational $\alpha$.

\smallskip

Furthermore, if $u\in S_{2k}(SL(2,\Z))$ as above, and $1\le \rho\le 2k-1$,
 applying to $\alpha$ the ``continued fractions trick'',
we can replace $(\alpha ,i\infty )$ by a sum of geodesic paths 
in the upper half--plane, joining pairwise cusps of the form
$g^{-1}(0)$ and $g^{-1}(i\infty )$, where $g$ varies in $SL(2,\Z )$, and then return
to $(0, i\infty )$ by transforming the integrand via
$\tau \mapsto g\tau$.  Thus, in particular, all values of (0.8),
corresponding to integer $\rho$'s inside the critical strip
and arbitrary characters $\kappa$, can be expressed as linear combinations
of modular symbols {\it with rational coefficients}, and
 span a finite--dimensional space over $\Q$.

\medskip

{\bf 0.4.  $p$--adic  Mellin--Mazur transform.}
Such expressions were used in [Ma2], [Ma3] in order
to produce a $p$--adic interpolation of values (0.8).
This problem  will make sense, if (after an appropriate normalization)
these values will lie in a finitely generated $\Z$--module,
so the  basic problem is to control the denominators.

\smallskip

As we already said, the main tool for such an interpolation was 
a  $p$--adic integral 
(the Mellin--Mazur transform) with respect to  a {\it $p$--adic 
pseudo--measure} (see below) constructed using modular symbols.
This transform  integrates  $\tau^{\rho -1}$ twisted by $\kappa$
against this pseudo--measure,
and for finite order $\kappa$
produces the classical values $L_{\kappa} (\rho )$
more or less by definition. (In fact, one works usually with
Dirichlet characters in place of $\kappa$, but the only difference
consists in the appearance of auxiliary Gauss sums).

\smallskip

Here are some details.

\smallskip

(a) {\it The $p$--adic integration domain and a naive pseudo--measure.} 
The following tentative construction applies to
any (absolutely convergent) series of the type (0.8) considered as a function
of variable $\kappa$ with fixed $\rho$. 

\smallskip

At the first approximation, consider $\Z_p$ with $\Z$
densely embedded  in it. The Boolean algebra of
closed/open subsets of $\Z_p$ is generated by
the {\it primitive subsets} $a+p^m\Z_p$, $m=0, 1, 2, \dots $;
$a\,\roman{mod}\,p^m$. Put 
$$
\mu_L(a+p^m\Z ):= \sum_{n\equiv a\,\roman{mod}\,p^m}
a_nn^{-\rho} .
\eqno(0.12)
$$ 
Any two primitive subsets either do not intersect, or one of them
is contained in the other. If one primitive subset $I$
is a disjoint union of a finite family of other primitive
subsets $I_j$, then $\mu_L(I)=\sum_j \mu_L(I_j).$
Thus $\mu_L$ extends to a $\C$--valued finitely additive function
on the Boolean algebra
of closed/open subsets of $\Z_p$. We will call such objects
{\it pseudo--measures} on $\Z_p$.

\smallskip

Generally, there is no chance that such a pseudo--measure
will tend $p$--adically to zero when $m\to\infty$,
even if its values lie in a finite--dimensional $\Q$--space.
As is explained in [Ma2], a Mazur's $p$--adic integral
of a function against such a pseudo--measure typically
 converges not because the smaller
primitive subsets have asymptotically vanishing pseudo--measure,
but because in a typical Riemann sum, {\it many approximately equal
 terms} of  not very large $p$--adic
size are involved, and  {\it the quantity
of summands $ \approx p^m$, tending to zero $p$--adically}, 
produces an unconventional non--Archimedean convergence effect. 

\smallskip

If the pseudo--measure of small subsets does
not tend to zero, the best one may hope for is that it will be
{\it bounded}, i.e. its values will lie in a
$\Z$--module of finite type. Even this usually will
not happen: for example, one can suspect that
$$
\mu_L(p^m\Z)=\sum_{n\equiv 0\,\roman{mod}\,p^m}
a_nn^{-\rho} = p^{-m\rho}\sum_n a_{np^m}n^{-\rho}
$$
will have denominator of order $p^{-m\rho}$.

\smallskip

A radical way to avoid this danger is to postulate
that $a_n=0$ if  $n$ is divisible by $p$. One can achieve
this cheaply, if $L$ admits an Euler product:
simply discard the $p$--th Euler factor of $L$.

\smallskip

(Notice an interesting Archimedean analogy:
the Mellin transform $\Lambda $ in (0.10)
produces $L$ {\it supplemented} by the initially
missing ``Euler factor at arithmetical infinity''.)

\smallskip

Returning to $L_{(p)}:=L$ {\it divided by its $p$--factor}, 
we may from now on look only at the group
of $p$--adic units $\Z_p^*\subset \Z_p$ by which
our pseudo--measure is now supported.

\smallskip

We repeat in conclusion, that the classical values (0.8) are tautologically
integrals of the locally constant function $\kappa$
against our pseudo--measure (0.12). (Of course, this is why chose it
in the first place). Only when we start
to interpolate and allow, say, continuous
$p$--adically valued  multiplicative characters in place of $\kappa$,
we will need the basics of such $p$--adic integration.

\smallskip 

{\it (b) Normalized  $p$--adic  pseudo--measure.} Let now $L$ be the Mellin transform
of an $SL(2,\Z)$--cusp form of weight $2k$ as above. Representing
the characteristic function of the set $a +p^m\Z$ by a linear
combination of the additive characters $\kappa$ modulo
$p^m$, and calculating $\Lambda_{\kappa}(\rho )$ as in 
(0.4), we see that $\mu_L(a+p^m\Z_p)$ is a linear combination of
modular symbols $\{bp^{-m}, i\infty \}$, $b\in \Z$. 

\smallskip
Conversely, we may
take an appropriate linear combination of such measures
and obtain the one that was used in [Ma2], [Ma3], namely
$$
\mu_p (a+p^m\Z_p) := \varepsilon^{-m}\{ap^{-m}, i\infty\}_k 
-p^{2k-2} \varepsilon^{-m+1} \{ap^{-m+1}, i\infty\}_k  .
\eqno(0.13)
$$
Here $\varepsilon$ is a root of  the (inverted) $p$--factor of $L$:
$\varepsilon^2-a_p\varepsilon +p^{2k-1}=0.$ If one of the two
roots is a $p$--adic unit, we get a bounded measure.
In any case, its growth can be controlled. The appearance of two summands
and $\varepsilon$ in (0.13) is a slightly more
sophisticated solution than the total discarding of the
$p$--th Euler factor.

\medskip

{\bf 0.6.  $\infty$--adic L\'evy--Mellin transform.} As it was suggested
in [MaMar2], let us make the following replacements in the picture
sketched above.

\smallskip

Replace $p$ by the arithmetic infinity  (archimedean 
valuation of $\Q$). Replace $\Z_p^*$ by the semi--interval $(0,1]$.

\smallskip

Call the classical Farey intervals with ends
$(g^{-1}(i\infty ), g^{-1}(0))$, $g\in SL(2,\Z)$, {\it primitive segments.}
They will be our replacement for the residue
classes $a+p^m\Z_p$. Exactly as residue classes, two open
primitive segments either do not intersect,
or one of them is contained in another. For an abelian group $W$,
call a pseudo--measure a $W$--valued finitely additive function
on segments with rational ends (see additional details below).

\smallskip

A typical pseudo--measure in this sense is the modular symbol
itself:
$$
\mu (\alpha ,\beta ) = \{\alpha ,\beta \}_k,
$$
in particular, $\mu (\alpha ,\infty ) = \{\alpha ,\infty \}_k$
which may be compared to (0.13).

\smallskip

As in the $p$--adic case, the pseudo--measure of a small segment
is not small in the archimedean sense. However, now we cannot hope
to compensate this by the non--Archimedean effect referred to above.

\smallskip

Instead, we suggest to use the following general feature
of our constructions:

\medskip

{\it (*) The Mellin transform of a cusp form, after  a  suitable 
normalization, can be naturally written as the  sum
over rational numbers in $(0,1]$ of values of a certain arithmetic function $a$:}
$$
A: =\sum_{\beta\in (0,1]\cap\Q} a(\beta ).
\eqno(0.14)
$$

\medskip

The values $a(p/q)$ involved here are essentially modular symbols
divided by a power of the denominator $q$. For details, see sec. 2 below.

\smallskip

Generally, a convergent series of the form (0.14) gives rise to an archimedean integral
in two related ways:

\smallskip
 
(i) {\it The first construction.}  We can define a pseudo--measure $\mu = \mu_a$  on the Boolean
algebra generated by segments with {\it irrational} ends in
$[0,1]$ putting  
$$
\mu (\alpha ,\beta ):= \sum_{\gamma \in (\alpha ,\beta )\cap\Q} a(\gamma ).
\eqno(0.15)
$$
so that  
$$
A=\int_0^1 d\mu .
\eqno(0.16)
$$
One can also treat (0.14) as a distribution on an appropriate space of test functions.

\smallskip

This is a direct analog of  (0.12),  however, it is not the version
that we will use in this paper.

\smallskip

(ii) {\it The second construction.} Let $r$ be a function defined on pairs of positive coprime 
integers $(p,q), p<q$ and decreasing sufficiently fast. For a real number $\xi$, 
denote by  $q_i(\xi )$  the denominator of
the $i$--th convergent to $\xi$, $i\ge 0.$
We can introduce  the L\'evy 1-form $l(\xi )d\xi$, associated to $r$,
and defined on $(0,\dfrac{1}{2}]$
by the following prescription:
$$
l(\xi )= l_r(\xi ) = \sum_{i=0}^{\infty} r (q_i(\xi ),q_{i+1}(\xi )) .
\eqno(0.17)
$$
According to a lemma by P.~L\'evy, for any pair $(p,q)$
as above, the set of all $\xi \in (0,\dfrac{1}{2}]$ for which there exists $i$
with $(p,q)=(q_i(\xi ),q_{i+1}(\xi ))$, fills a primitive semi--interval 
of length $\dfrac{1}{(p+q)q}$. Moreover, this $i$ is uniquely defined.
Therefore, when $r(p,q)$ decreases sufficiently rapidly to assure
convergence, we get
$$
\int_0^{1/2} l_r(\xi )d\xi =\sum_{\alpha=p/q\in (0,1]} \frac{r(p,q)}{(p+q)q}.
\eqno(0.18)
$$
In particular, we get $A$ from (0.14) if we choose
$$
r(p,q) := a(p/q)(p+q)q.
\eqno(0.19)
$$
When $A$ comes from a modular form (classical or Maass),
so that the summands  $a(\beta )$ are concocted  of (classical or Maass)
modular symbols, we will call the integral in (0.12)
{\it the $\infty$--adic L\'evy--Mellin transform.}

\smallskip

The L\'evy functions and their generalizations appear
also in a different context: that of linearizations of the germs of
analytic diffeomorphisms of one complex variable $z$ with an indifferent
fixed point. For example, a germ with linear part 
$ e^{2\pi i\xi}z$ is linearizable iff the Brjuno number of $\xi$
$$ 
b(\xi ):= \sum_{i=0}^{\infty} \frac{\roman{log}\,q_{n+1}}{q_n}
$$
is finite . In fact,
an interesting theory is developed/reviewed in the papers
[MarMouYo1] and [MarMouYo2] for another Brjuno
function $B(\xi )$, which differs from $b(\xi )$ by $O(1)$,
but satisfies a functional equation and has a complex version
closely resembling some constructions
in the theory of modular forms.
In our context, it can be used
for calculation of the derivative of some classical $L$--series
at certain points. This looks like an interesting variation on the subject
of L\'evy--Mellin transform.

\medskip

{\bf 0.7. A summary of $p$--adic/$\infty$--adic analogies.} For clarity,
we summarize the suggested analogies in the following lines:

\medskip

$\phantom{0000000000000000000}\Z_p^*  \phantom{000000} \Longleftrightarrow  \phantom{0000000} (0,1] $

\smallskip

$\phantom{0000000000000000000}\cup \phantom{00000000000000000000} \cup$

\smallskip

$\phantom{0000000000000000000}\Z\phantom{0000000} \Longleftrightarrow  \phantom{0000000}\Q\cap (0,1]$

\medskip

$\phantom{00000000000000} a+p^m\Z_p \phantom{00000} \Longleftrightarrow  \phantom{0000000} primitive\  (Farey) \ segments$

\medskip

$\phantom{00000000000000}   \sum_{m=1}^{\infty}\dfrac{a_m}{m^{\rho}}    \phantom{0000} \Longleftrightarrow \phantom{000000} \sum_{0<p/q\le 1}\dfrac{a(p/q)}{q^{\rho}} \phantom{0000000}  $

\medskip

{\it Mazur--Mellin transform}$ \phantom{000000} \Longleftrightarrow  \phantom{0000000}$
{\it L\'evy--Mellin\ transform}

\bigskip

\centerline{ \bf \S 1.  Pseudo--measures}

\smallskip

\centerline{\bf  associated with period--like functions}

\medskip

{\bf 1.1. An heuristic construction.} For the moment, we adopt the viewpoint
of [LZ2], Ch. II, sec. 5. Fix a complex number $s$ such that $s(1-s)$ is an eigenvalue
of the standard hyperbolic Laplace operator on $\bold{C}$
producing a $PSL(2,\bold{Z})$--invariant Maass wave form $u(z)=u_s(z), z\in H$.
Define complex powers by the usual formula
$t^s := e^{s\roman{log}\,t}$ where the branch of the logarithm is determined by
the normalization $-\pi <\roman{arg}\,t\le \pi.$
As is shown in [LZ1], there exists a distribution $U(t)=U_s(t)$ on $\bold{R}$,
whose  values on the following  test functions of $t$ ($z$ enters as a parameter)
$$
(\roman{Im}\, z)^s |z-t|^{-2s},\quad (z-t)^{-2s},\quad \chi_{(-\infty ,0)}(t) (z-t)^{-2s}
$$
are respectively $u(z)$ (the initial Maass form), a function $f(z)$ holomorphic in $\bold{C}\setminus
\bold{R}$, and a {\it period function} $\psi (z)$ defined and holomorphic
in $\bold{C}^{\prime}:=\bold{C}\setminus (-\infty ,0].$ Here $\chi$ is the characteristic function of $\bold{R}_{-}$,
in other words
$$
\psi (z)=\int_{-\infty}^0 (z-t)^{-2s} U(t)dt .
\eqno(1.1)
$$
The distribution $U$ is automorphic in the following sense:
for all 
$$
g= \left( \matrix
a & b\\
c &  d
\endmatrix \right) \in SL(2,\bold{Z})
$$
we have
$$
U\left(\frac{at+b}{ct+d}\right)= |ct+d|^{2-2s}  U(t).
\eqno(1.2)
$$
Thus, (1.1) has  the same structure as (0.1).

\smallskip

Consider now only $g\in SL(2,\bold{Z})$ {\it with non--negative entries}.
Then for any $z\in \bold{C}^{\prime}$ we have  also $gz:=\dfrac{az+b}{cz+d}\in \bold{C}^{\prime}$.
From (1.2) we find formally
$$
\psi (gz)=\int_{-\infty}^0(gz-t)^{-2s}U(t) =\int_{g^{-1}(-\infty)}^{g^{-1}(0)}
(gz-g\tau)^{-2s}U(g\tau )\,d(g\tau ).
\eqno(1.3)
$$
The direct calculation using (1.2) reduces the integrand  to the form
$$
\left[\frac{z-\tau}{(cz+d)(c\tau +d)}\right]^{-2s} |c\tau +d|^{-2s+2}U(\tau )\frac{d\tau}{|c\tau 
+d|^2}.
\eqno(1.4)
$$
Since $a\ne 0$, we have
$$
g^{-1}(-\infty )= -\frac{d}{c}< -\frac{b}{a}=g^{-1}(0),
$$
and  hence for $\tau \in (g^{-1}(-\infty ), g^{-1}(0))$
we  have $c\tau +d > 0$.  This shows that all terms involving $c\tau+d$ in (1.4)
cancel, so that finally we find formally
$$
\psi (gz) = (cz+d)^{2s} \int_{-d/c}^{-b/a} (z-\tau)^{-2s}U(\tau )d\tau .
\eqno(1.5)
$$
Thus if $(\alpha ,\beta ) =(g^{-1}(-\infty ),g^{-1}(0))$ with $g$ as above, and if we put
$$
\mu (\alpha ,\beta ) (z):= (cz+d)^{-2s} \psi (gz) =\int_{\alpha}^{\beta}(z-t)^{-2s}U(t)dt,
\eqno(1.6)
$$
then for three intervals of this type $(\alpha,\beta )$, $(\beta ,\gamma )$, $(\alpha ,\gamma )$
we would have from (1.6)
$$
\mu (\alpha ,\beta ) (z) +\mu (\beta ,\gamma )(z) =
\mu (\alpha ,\gamma ) (z).
\eqno(1.7)
 $$
 As we will see, all primitive intervals in $\bold{R}_{-}$ are of this form,
 so that we have formally constructed a pre--measure (see below) on the (left half of)
 $\bold{P}^1(\R)$,  extendable  to a pseudo--measure on this half
 with values in the space of holomorphic functions
 on $\bold{C}^{\prime}$,
 in view of [MaMar2], Theorem 1.8.
 
 \smallskip
 
 The weak point of this reasoning, about which the word
 ``formally'' is supposed to warn the reader, is this:
the  functions $\chi_{(\alpha ,\beta)} (t) (z-t)^{-2s}$ generally do not 
belong to the space of test functions as it is defined in [LZ1], p. 225.
Therefore the integrals  in the r.h.s. of (1.5), (1.6) {\it a priori} make no sense.

\smallskip

Our heuristic reasoning is in fact  a simple extension of the
formal argument on p. 222 of [LZ2], ``proving'' the three--term functional 
equation for $\psi (z)$. 

\smallskip

In the next subsections, we will provide a precise construction
of the pseudo--measures, whose
values on the intervals considered above are given by
$$
\mu (g^{-1}(-\infty ),g^{-1}(0))(z):= (cz+d)^{-2s}\psi (gz)
\eqno(1.8)
$$
without appealing to the integral representation (1.6), but making use of the theory developed 
in [LZ1].

\medskip

{\bf 1.2. Preliminaries I: left primitive segments.} As in [MaMar2],   
we consider $\Q\subset \R\subset \C$
as points of an affine line with a fixed coordinate $z$.
Completing this line by one point $\infty=-\infty=i\infty$, we get
points of the projective line $\P^1(\Q )\subset \P^1(\R )\subset \P^1(\C )$
(Riemannian sphere). $GL(2,\bold{C})$ acts on  $\P^1(\C )$
by fractional linear transformations.
{\it Segments}  are defined as non--empty 
connected subsets of $\P^1(\R )$. A segment is called infinite if $\infty$ is in its closure,
otherwise it is called finite. The boundary of each segment generally
consists of an unordered pair of points $(\alpha ,\beta )$ in $\P^1(\R )$.
We will identify a segment with an {\it ordered} pair of its ends:
the additional element of structure is its orientation from $\alpha$
to $\beta$. For our purposes, it is usually inessential
whether one or two boundary points belong
to the segment. In this section we will consider mostly
{\it left} segments, that is, ones for which $-\infty\le\alpha ,\beta \le 0.$
One--point segments are sometimes  called improper ones.

\smallskip

 A  segment
is called {\it rational} if its ends are in $\P^1(\Q )$, and {\it primitive}, or Farey,
if it is of the form $(g(\infty ),g(0))$ for some $g\in GL(2,\bold{Z})$.
 
\smallskip 

In [MaMar2] we called {\it a pseudo--measure with values in an abelian group $W$}
a finitely additive $W$--valued function on the Boolean algebra
of rational segments, vanishing on improper segments. We extended it
to oriented segments by the condition that $\mu (\alpha ,\beta )=
-\mu (\beta ,\alpha ).$ 

\smallskip

In this section, we will construct pseudo--measures supported by left segments.
As it follows from [MaMar2], each such pseudo--measure
is defined by its restriction to the set $P$ of positively oriented left primitive segments.
We will use the following enumeration of the latter.

\smallskip

Denote by $S\subset SL(2,\bold{Z})$ the sub--semigroup of matrices with
non--negative entries $a,b,c,d.$ For any $g\in S$, 
$(g^{-1}(-\infty ),g^{-1}(0))$ is in $P$. In fact, if $c\ne 0$,
$$
g^{-1}(\infty )= \frac{d}{-c} < g^{-1}(0) =\frac {-b}{a} ,
$$ 
because $ad-bc=1$. If $c=0$, then $a=d=1$, and again
$$
g^{-1}(\infty )= -\infty < g^{-1}(0) =-b.
$$ 
Finally, the case $a=0$ does not occur in $S$.

\smallskip

One easily sees that this map $S\to P$: $g\mapsto (g^{-1}(-\infty ),g^{-1}(0))$
is in fact a bijection.

\medskip

{\bf 1.3. Preliminaries II: the slash operators of complex weight.}
Here we summarize the considerations of [LZ2], p. 240, and  [HiMaMo], sec 3.
They determine a partial map
$$
(\varphi , g) \mapsto \varphi |_{s}g
\eqno(1.9)
$$
allowing us to make sense of and correctly calculate expressions as those
appearing in (1.4), (1.6). 

\smallskip

For proofs, see [HiMaMo].

\medskip

{\it (i) Definition domain.} The argument $\varphi =\varphi (z)$ in  (1.9)
can be an arbitrary function holomorphic in some domain
of the form  $\bold{C}\setminus (-\infty ,r]$,  $r\in \bold{R}$.
Such functions form a $\bold{C}$--algebra which we will denote $\Cal{F}$.
Period functions $\psi = \psi_s$ belong to  $\Cal{F}.$

\smallskip

In [HiHaMo], any point $r$ such that $\varphi \in \Cal{F}$ is holomorphic
in  $\bold{C}\setminus (-\infty ,r]$, is called {\it a branching point}
of $\varphi$.

\smallskip

The argument $g$ in (1.9) can be an arbitrary $(2,2)$--matrix $g$ with integer
entries $(a,b,c,d)$ and non--zero determinant such that either $c>0$, or $c=0; a,d>0.$
Denote by $\Cal{G}$ the set of such matrices.
The set $S$ describing left primitive segments in 1.2 is a subset of $\Cal{G}$.
When $g\in \Cal{G}$ and $s\in \bold{C}$, the function $(cz+d)^s$
belongs to $\Cal{F}.$

\smallskip

A pair $(\varphi ,g )\in \Cal{F} \times \Cal{G}$ belongs to the definition
domain $\Cal{D}\Cal{S}$ of the slash operator (1.9) if $\varphi$
admits a branching point $r$ such that either
$a-cr >0$, or $a-cr=0$ and $dr-b<0.$  For a period function $\phi =\psi$,
we can take $r=0$, and $g$ will do if $a>0$ or else $a=0, b>0.$ 

\smallskip

Let $\Cal{G}^+$ be the set of matrices in $\Cal{G}$ such that $b,d\ge 0$
and either $a>0$, or $a=0, b>0.$ Again, $S\subset \Cal{G}^+$. Denote by
$\Cal{F}_0$ the subspace of $\Cal{F}$ admitting $0$ as a branch point.
Then $\Cal{F}_0\times \Cal{G}^+ \subset \Cal{D}\Cal{S}.$

\medskip

{\it (ii) Slash operator of weight $s$.} It is the map $\Cal{D}\Cal{S}\to \Cal{F}$
defined by 
$$
(\varphi (z) , g) \mapsto (\varphi |_{s}g)(z) := |\roman{det}\,g|^s (cz+d)^{-2s}
\varphi (gz).
\eqno(1.10)
$$
It is well defined. Moreover,
it sends $\Cal{F}_0\times \Cal{G}^+$ to $\Cal{F}_0.$

\smallskip

{\it (iii) Properties of the slash operator.} The basic property is that slash
operator is a honest action: if $g_1, g_2 \in \Cal{G}$ and $(\varphi ,g_1),
(\varphi |_sg_1,g_2), (\varphi ,g_1g_2)\in \Cal{D}\Cal{S}$, then
$$
\varphi |_s(g_1g_2)=( \varphi |_sg_1)|_sg_2.
$$
(Formally, it is the associativity of the triple product of $\varphi, g_1,g_2.$)
Applying this to $\Cal{F}_0\times \Cal{G}^+$, one can check that $|_s$ defines a
right action
of the multiplicative semigroup $\Cal{G}^+$  on $\Cal{F}_0$ ([HiMaMo],
Remark 3.4).

\smallskip

From 1.2 one sees that if $(g^{-1}(-\infty ),g^{-1}(0))$ is a left primitive  segment,
then $g\in \Cal{G}^+$. Since $\psi \in \Cal{F}_0$ in the first equality
of (1.6), this expression for $\mu (\alpha ,\beta )(z)$ (disregarding the second equality and the poorly
defined  integral) makes sense,
and the slash action can be further iterated.

\medskip

{\bf 1.4. The pre--measures related to  period--like functions.} 
Choose a complex number $s$ and a function $\psi (z)\in \Cal{F}_0$ satisfying 
the three term functional equation
$$
\psi (z) =\psi (z+1)+ (z+1)^{-2s}\psi \left(\frac{z}{z+1}\right)
\eqno(1.11)
$$
Thus, $\psi$ is a period--like function in the sense of [LZ], Ch. III.

\smallskip

For a left primitive segment $(\alpha ,\beta )=(g^{-1}(-\infty ),g^{-1}(0))$, put 
 $$
\widetilde{ \mu} (\alpha ,\beta)(z) = (cz+d)^{-2s} \psi (gz) =\psi|_s(z).
\eqno(1.12)
$$

Consider now three left primitive
segments $(\alpha, \beta ) = (g_1^{-1}(-\infty ),g_1^{-1}(0))$, 
$( \beta ,\gamma ) = (g_2^{-1}(-\infty ),g_2^{-1}(0))$, $(\alpha, \gamma ) = (g_3^{-1}(-\infty ),g_3^{-1}(0))$.
In plain words, the third segment is broken into two others by
a point $\beta$ in the middle.
 
 \medskip
 
 {\bf 1.4.1.  Lemma.}  {\it We have
 $$
\widetilde{\mu} (\alpha ,\beta ) (z) +\widetilde{\mu} (\beta ,\gamma )(z) =
\widetilde{\mu} (\alpha ,\gamma ) (z).
\eqno(1.13)
 $$}

{\bf Proof.}  {\it Case 1: $(\alpha ,\beta ,\gamma ) = (-\infty , -1, 0)$}. In this case 
$$
g_1 = T:= \left( \matrix 1 & 1\\ 0 &  1 \endmatrix \right), \quad
g_2 = T^{\prime}:= \left( \matrix 1 & 0\\ 1&  1 \endmatrix \right), \quad
g_3 = I:= \left( \matrix 1 & 0\\ 0 &  1 \endmatrix \right) ,
$$
and the equation (1.13) coincides with (1.11) which can be written as
$$
\psi |_sI =\psi |_s T +\psi|_s T^{\prime} .
\eqno(1.14)
$$

\smallskip

{\it Case 2:  $g_1=Tg, g_2=T^{\prime}g, g_3=g$}, where $g\in SL(2,\bold{Z})$ 
is a matrix with non--negative entries. 

\smallskip

In this case, (1.13) reads
$$
\psi |_sg =\psi |_s Tg +\psi|_s T^{\prime}g
$$
which obviously holds in view of  (1.14) and the associativity of the slash operator
restricted to $\Cal{F}_0\times \Cal{G^+}$.

\smallskip

{\it General case.} In fact, the previous case {\it is} general: 
we  necessarily have $g_1= Tg_3$ and $g_2=T^{\prime}g_3.$

\smallskip

Let us check this for the case when $\alpha \ne -\infty$ leaving 
the remaining case to the reader. Put $\alpha = \dfrac{d}{-c},
\gamma = \dfrac{-b}{a}$ as in 1.2 where $(a,b, c, d)$ are the entries of $g_3$.
Then the only possible value of $\beta$ is
$\beta = \dfrac{-(b+d)}{a+c}=\dfrac{b+d}{-(a+c)}$ as is well known from
the classical theory of Farey series. This fact directly translates into
$g_1=Tg_3, g_2=T^{\prime}g_3.$

\smallskip

This completes the proof of the lemma.

\smallskip

{\it Remark.} Notice that if $\psi (z)$ is an actual period function for a Maass wave form, 
the lemma becomes obvious in view of the integral representation
of $\psi (z)$ proven in [LZ], Ch. II, sec. 1. The relevant formula on the p. 212
of [LZ] reads (we have replaced the notation $\psi_1(\zeta )$ by $\psi (\zeta )$
and changed the sign):
$$
(c\zeta +d)^{-2s} \psi (g\zeta ) = \int^{g^{-1}(0)}_{g^{-1}(\infty )}\{u, R^s_{\zeta}\}(z).
\eqno(1.15)
$$
In this formula, we integrate a closed form along an arbitrary path leaving 
$\zeta$ and $\overline{\zeta}$ to the right of it. Additivity (1.13)
becomes evident.

\smallskip

We will use this integral representation in the next section.

\medskip

{\bf 1.4.2. The pre--measure on left segments.} To define a pre--measure 
in the sense of [MaMar], {\it supported by the subset of left primitive segments},
it remains to complete the definition (1.12) of the function
$\widetilde{\mu}$ by putting for $\alpha <\beta \le 0$
$$
\widetilde{\mu}(\beta ,\alpha ):= -\widetilde{\mu}(\alpha, \beta ), \quad \widetilde{\mu}(\alpha ,\alpha)=0.
$$
\medskip
One easily checks that (1.13) continues to hold on this extended domain.

\medskip

{\bf 1.5. The pseudo--measure related to a period--like function.}  Now we can state the 
main result of this section.

\medskip

{\bf 1.5.1. Theorem.} {\it There exists a unique finitely additive
function $\mu$ with values in $\Cal{F}_0$ coinciding with $\widetilde{\mu}$
on left primitive segments and vanishing on all rational segments in $(0,\infty ).$}

\smallskip

{\bf Sketch of proof.} We simply recall the plan of proof of the 
Theorem 1.8 of [MaMar2].  It  consists of the following steps. 

\smallskip

1) Using the ``continued fractions trick'',  we show that for any non--positive rational (or infinite)
 $\alpha ,\beta$ one can find a sequence of rational non--positive numbers
 $\alpha_0=\alpha, \alpha_1,\dots , \alpha_n=\beta$ such that $(\alpha_i,\alpha_{i+1})$
 is a left primitive segment for all $i=0, \dots , n-1.$ Such sequence is called
 a primitive chain connecting $\alpha$ to $\beta$.
 
 \smallskip
 
 2) Having chosen such a primitive chain, we put 
 $$
 \mu (\alpha ,\beta ):= \sum_{i=0}^{n-1} \widetilde{\mu} (\alpha_i,\alpha_{i+1}).
 \eqno(1.16)
 $$
 
 \smallskip
 
 3) The fact that (1.16) does not depend on the choice of the connecting primitive
 chain is checked by proving that any two chains can be transformed one to another
 by using ``elementary moves'' compatible with relations that hold for $\widetilde{\mu}$.
 An elementary move essentially replaces a Farey interval
 $\left(\dfrac{a}{c},\dfrac{b}{d}\right)$ by the chain $\left(\dfrac{a}{c},\dfrac{a+b}{c+d}\right)$,
 $\left(\dfrac{a+b}{c+d},\dfrac{b}{d}\right)$,
 or vice versa.

 \smallskip
 
 4) Finally, we have to check that (1.16) implies finite additivity
 and the sign change  after the change of orientation. This is straightforward.
 
 \bigskip
 
 {\bf 1.6. Modularity.} Let $\Gamma$ be a subgroup of $SL(2,\Z)$,
 $W$ a left $\Gamma$--module.
 
 \smallskip
 
 In [MaMar2], a pseudo--measure $\mu$ with values
 in $W$ is called {\it $\Gamma$--modular}, if for all $g\in \Gamma$
 and $\alpha ,\beta \in \bold{P}^1(\bold{Q})$ we have
 $$
 \mu (g\alpha ,g\beta ) =g\mu (\alpha ,\beta ).
 $$
 It was checked that such pseudo--measures correspond  to parabolic 1--cocycles.
 
 \smallskip
 In our context, this is replaced by the following property:
 for all $g$ with  $g^{-1} \in S$
 and any left segment $(\alpha ,\beta )$,
 $$
 \mu (g^{-1}(\alpha ),g^{-1}(\beta )) =\mu (\alpha ,\beta )|_sg.
 \eqno(1.17)
 $$
 In fact, it suffices to check this for left primitive segments $(\alpha ,\beta )=
 (h^{-1}(-\infty ),h^{-1}(0))$, in which case we have
 $$
  \mu (g^{-1}(\alpha ),g^{-1} (\beta )) =  \mu ((hg)^{-1}(-\infty ),(hg)^{-1} (0 ))=
  $$
  $$
 \psi |_s( hg)=(\psi |_sh)|_sg = \mu (\alpha ,\beta )|_sg.
 $$
 Since the right slash action of $g$ can be considered as the left
 action of $g^{-1}$, we can say that (1.17) expresses the modularity of $\mu$ with
 respect to the multiplicative semigroup $S^{-1}\subset SL(2,\Z ).$

\bigskip

\centerline{\bf \S 2. Maass $L$--functions}

\smallskip

\centerline{\bf and their Mellin--L\'evy transforms}

\medskip

{\bf 2.1. Maass $L$--series as sums over rational numbers.}
Let $u=u_s$ be a Maass cusp form, which is an eigenfunction
with respect to all Hecke operators
$$
T_m := \sum_{ad=m \atop 0< b\le d} 
\left( \matrix
a &- b\\
0 &  d
\endmatrix \right)
\eqno (2.1)
$$
acting via slash operator of weight 0: $u\mapsto u|_0T_m =\lambda_mu.$ 
\smallskip

Put
$$
L_u(\rho ):= \sum_{m=1}^{\infty}\frac{\lambda_m}{m^{\rho}}.
\eqno(2.2)
$$
The action of the Hecke operators on $u$ induces an action
on the period functions, which can be explicitly described by a nice
formula,
for example, as in [M\"uh]. However, we will need a different expression,
involving the pseudo--measure $\mu_u$, and we will start with
an heuristic derivation of it, as in 1.1.

\smallskip

Let us formally apply the slash operator  $|_{-s}$ (see (1.10))
to the boundary measure $U(t)dt$ and denote the resulting action upon
the respective period function $\psi$ by $T^*_m$. In this heuristic
calculation we ``define'' $\psi$ by (1.1). The choice of weight $-s$
is motivated by the invariance property (1.2). We get:
$$
(\psi|T^*_m )(\zeta):=\int_{-\infty}^0 (\zeta-t)^{-2s} (U(t)\,dt\,|_{-s}T_m) =
$$
$$
 \sum_{ad=m\atop 0< b\le d}\left(\frac{d}{a}\right)^{s}\int_{-\infty}^0 (\zeta-t)^{-2s} U\left(\frac{at-b}{d}\right)\,
d\left(\frac{at-b}{d}\right) .
$$
Make the change of variable $\tau =\dfrac{at-b}{d}$. The last integral takes form
$$
 \sum_{ad=m\atop 0< b\le d} \left(\frac{d}{a}\right)^{s}\int_{-\infty}^{-b/d} \left(\zeta-\frac{d\tau+b}{a}\right)^{-2s} U(\tau )\,d\tau=
$$
$$
\sum_{ad=m\atop 0< b\le d}\left(\frac{d}{a}\right)^{s}\int_{-\infty}^{-b/d} \left(\frac{dz+b}{a}-\frac{d\tau+b}{a}\right)^{-2s} U(\tau )\,d\tau ,
$$
where $z =\dfrac{a\zeta-b}{d}.$ The integral in the last sum can be rewritten as
$$
\left( \frac{a}{d}\right)^{2s} \int_{-\infty}^{-b/d} (\zeta -\tau )^{-2s} U(\tau )\,d\tau .
$$
Thus heuristically
$$
(\psi|T^*_m)(\zeta )=( \mu (-\infty, 0)|T^*_m)(\zeta ) = \sum_{ad=m\atop 0< b\le d} \left(\frac{a}{d}\right)^s
 \mu \left(-\infty,-\frac{b}{d}\right)\left(\frac{a\zeta -b}{d}\right) =
$$
$$
\sum_{ad=m\atop 0< b\le d}
\mu\left. \left(-\infty,-\frac{b}{d}\right)\right|_s \left( \matrix
a &- b\\
0 &  d
\endmatrix \right)(\zeta) .
\eqno(2.3)
$$
This expression is useful for our purposes because it 
allows us to represent (the somewhat normalized)
Dirichlet series $L_u(s)$ as a natural sum
over rational numbers. We will state now the respective theorem:

\medskip

{\bf 2.2. Theorem.} {\it We have }
$$
\psi (z) \sum_{m=1}^{\infty} \frac{\lambda_m}{m^{\rho}}=
\zeta (\rho -s) \zeta(\rho +s) \sum_{q=1}^{\infty} \frac{1}{q^{\rho}}
\sum_{0\le p<q\atop  (p,q)=1} \mu (-\infty ,-p/q)
|_s\left(\matrix  1&-p\\0&q  \endmatrix \right) (z) .
\eqno(2.4)
$$
\smallskip

{\bf Proof.} {\it Step 1.} First, we have to supply an honest proof of 
(2.3). In [LZ2], Ch. II, sec. 2, the authors construct
a differential 1--form $\{u, R^s_\zeta\} (z)$  which we invoked at the end of 1.4.1. 
It has the following properties:

\medskip

(i) $\{u, R^s_\zeta\} (z)$  is a closed smooth form of $z$ varying in the complex
upper half--plane $H$.
It depends on the parameter $\zeta\in\C$ holomorphically when $z\ne \zeta,
\overline{\zeta}$.
Generally, it is multivalued, but a well defined  branch can be chosen
on the complement in $H$ of a path joining $\zeta$ to $\overline{\zeta}.$

\smallskip

(ii) The period function $\psi (\zeta),\zeta \in H$ for $u$ (up to a constant proportionality factor) can
be then written as an integral
$$
\psi (\zeta ) =\int^0_{-\infty} \{u, R^s_\zeta\} (z)
\eqno(2.5)
$$
taken along any path in $H$ leaving $\zeta$ to the left of it.

\smallskip

Let us now assume that $u|_0T_m=\lambda_mu$ for $T_m$ from (2.1) and a constant $\lambda_m$.
Then we have from (2.5) and (2.1)
$$
\lambda_m\psi (\zeta ) = \int^0_{-\infty} \left\{ \sum_{ad=m \atop 0< b\le d} 
u\left( \frac{az-b}{d}\right), R^s_\zeta\right\} (z) .
\eqno(2.6)
$$
For each $a,b,d$ fixed, we first want to make the implicit argument $z$ of $R_{\zeta}^s$
the same as that of $u$, that is, $\dfrac{az-b}{d}$.  We have (see [LZ2], p. 211):
$$
R_{\zeta}(z)=\frac{i}{2}((z-\zeta )^{-1} -(\overline{z}-\zeta )^{-1}) =
$$
$$
\frac{a}{d}\cdot \frac{i}{2}\left(\left(\frac{az-b}{d}-\frac{a\zeta -b}{d}\right)^{-1} - \left(\frac{a\overline{z}-b}{d}-
\frac{a\zeta -b}{d}\right)^{-1}\right)=
\frac{a}{d} R_{\xi}\left(\frac{az-b}{d}\right), 
$$
where
$\xi : = \frac{a\zeta -b}{d}$ . 

\smallskip
Substituting this into (2.6), we obtain:
$$
\lambda_m\psi (\zeta ) =\sum_{ad=m \atop 0< b\le d}\left(\frac{a}{d}\right)^s \int^0_{-\infty} \left\{ 
u\left( \frac{az-b}{d}\right), R^s_{\frac{a\zeta -b}{d}}\left(\frac{az-b}{d}\right)\right\}  .
\eqno(2.7)
$$

Considering now $z\mapsto \frac{az-b}{d}$ as a holomorphic change of variables,
we infer from the Lemma on p. 210 of [LZ2] that the integrand in  the respective term
of (2.7) can be rewritten as
$$
\{u, R^s_{\frac{a\zeta -b}{d}}\}\left(\frac{az-b}{d}\right) .
$$
Hence finally
$$
\lambda_m\psi (\zeta ) =\sum_{ad=m \atop 0\le b<d}\left(\frac{a}{d}\right)^s \int^{-b/d}_{-\infty} 
\{u, R^s_{\frac{a\zeta -b}{d}}\} (z) 
=
$$
$$
 \sum_{ad=m \atop 0< b\le d} \mu (-\infty ,-b/d)|_s \left(
\matrix  a&-b\\0&d  \endmatrix \right) (\zeta ) .
\eqno(2.8)
$$
This is formula (2.3), written for $u$ which is an eigenfunction of $T_m$,
and its respective period function.

\medskip

{\it Step 2.} Multiply now the identity (2.8)
by $m^{-\rho}$ and sum over all $m=1,2, \dots $ . Again replacing the 
notation of the free variable $\zeta$ by $z$ (in order not to confound it with
Riemann's zeta in (2.11) below), we obtain 
$$
\psi (z) \sum_{m=1}^{\infty} \frac{\lambda_m}{m^{\rho}} =\sum_{m=1}^{\infty}
\frac{1}{m^{\rho}}
\sum_{ad=m \atop 0< b\le d} \mu (-\infty ,-b/d)
|_s\left(\matrix  a&-b\\0&d  \endmatrix \right) (z ) .
\eqno(2.9)
$$
Each matrix in (2.9) can be uniquely written in the following way:
$$
\left(\matrix  a&-b\\0&d  \endmatrix \right) =
\left(\matrix  d_2&-pd_1\\0&qd_1  \endmatrix\right)=
\left(\matrix  1&-p\\0&q  \endmatrix \right)
\left(\matrix  1&0\\0&d_1  \endmatrix \right) 
\left(\matrix  d_2&0\\0&1  \endmatrix \right) \, ,
\eqno (2.10)
$$
where $m=d_1d_2q$, $d_i\ge 1$, $0< p\le q$, $(p,q)=1.$
Moreover, arbitrary quadruple $(d_1,d_2,p,q)$ satisfying
these conditions produces one term in (2.9).

\smallskip

From (2.10) and the associativity of the slash operator (1.10) it follows that
$$
|_s\left(\matrix  a&-b\\0&d  \endmatrix \right)= 
|_s\left(\matrix  1&-p\\0&q \endmatrix \right) \cdot  d_1^{-s}d_2^{s}.
$$
Hence we can rewrite (2.9) as follows: 
$$
\psi (z) \sum_{m=1}^{\infty} \frac{\lambda_m}{m^{\rho}} =\sum_{q,d_1,d_2=1}^{\infty}
\frac{1}{q^{\rho}d_1^{\rho -s}d_2^{\rho+s}}
\sum_{0< p\le q\atop  (p,q)=1} \mu (-\infty ,-p/q)
|_s\left(\matrix  1&-p\\0&q  \endmatrix \right) (z ) =
$$
$$
\zeta (\rho -s) \zeta(\rho +s) \sum_{q=1}^{\infty} \frac{1}{q^{\rho}}
\sum_{0< p\le q\atop  (p,q)=1} \mu (-\infty ,-p/q)
|_s\left(\matrix  1&-p\\0&q  \endmatrix \right) (z) .
\eqno(2.11)
$$

{\bf 2.3. L\'evy--Mellin transform.} Put now
$$
r_u(p,q):= (p+q)q^{1-\rho} \mu (-\infty ,-p/q)
|_s\left(\matrix  1&-p\\0&q  \endmatrix \right) (z)\cdot \psi (z)^{-1} .
$$
and
$$
l_u(\xi ):= \sum_{i=0}^{\infty} \sum_{i=0}^{\infty} r (q_i(\xi ),q_{i+1}(\xi )) .
$$
From (2.11) and (0.18) we get the following

\medskip 

{\bf 2.3.1. Corollary.} {\it Let $u$ be a Maass cusp form, $\Delta u= s(1-s)u$,
$u|T_m=\lambda_mu$ for all $m\ge 1.$ Put
$$
\Lambda_u(\rho ):= \zeta (\rho -s)^{-1}\zeta (\rho +s)^{-1} \sum_{m=1}^{\infty}
\frac{\lambda_m}{m^{\rho}}.
$$
Then}
$$
\Lambda_u(\rho ) = \int_0^{1/2} l_u(\xi )d\xi.
$$

\smallskip

{\bf 2.3.2. Remark.} The class of series of the form (0.18) involving modular symbols
includes also D.~Goldfield's Eisenstein series, cf. [Go2]. They certainly deserve
further study.

\bigskip

{\bf 2.4. Hecke operators on period functions via continued fractions.}
Consider the sequence of normalized convergents $b/d$
as in [MaMar2], (1.5). When $0< b/d<1$, it starts with
$$
-\infty=\frac{1}{0}=:\frac{b_{-1}}{d_{-1}},\
0=\frac{0}{1} = :\frac{b_0}{d_0},\ \dots ,\  b/d=\frac{b_n}{d_n},
$$
where $n=n(b/d)$ is the length of the continued fraction expansion.
\bigskip

The following sequence of left primitive segments
$I_k=I_k(b/d)$ connects $-\infty$ to $-b/d$.
We order their ends from the left one to the right one, and 
put minus before those that should be run in the opposite direction 
in our chain:
$$
I_0=(-\infty ,0)=\left(-\frac{b_{-1}}{d_{-1}},- \frac{b_0}{d_0}\right),\
I_1=-\left(-\frac{b_{1}}{d_{1}},- \frac{b_0}{d_0}\right),\
$$ 
$$
I_2=\left(-\frac{b_{1}}{d_{1}},- \frac{b_2}{d_2}\right),\ 
I_3=-\left(-\frac{b_{3}}{d_{3}},- \frac{b_2}{d_2}\right),\
$$
and generally
$$
I_k=(-1)^k\left(-\frac{b_{k-\varepsilon_k}}{d_{k-\varepsilon_k}},
-\frac{b_{k-\varepsilon_{k+1}}}{d_{k-\varepsilon_{k+1}}}
  \right)
$$
where $\varepsilon_{k}=1$ for even $k$ and $0$ for odd $k$.
\smallskip

This means that 
$$
(-1)^kI_k=(g_k^{-1}(-\infty ),g_k^{-1}(0))
\eqno(2.12)
$$
where
$$
g_k=g_{k,b/d} = \left( \matrix
d_{k-\varepsilon_{k+1}} & b_{k-\varepsilon_{k+1}} \\
d_{k-\varepsilon_{k}} &  b_{k-\varepsilon_{k}}
\endmatrix \right) \in S\,.
\eqno(2.13)
$$
Therefore, (2.8) can be rewritten as
$$
\lambda_m\psi (\zeta ) =\sum_{ad=m \atop 0< b\le d}\left(\frac{a}{d}\right)^s
\sum_{k=0}^{n(b/d)} (-1)^k \int^{0}_{-\infty} 
\{u(g_{k,b/d}(z)), R^s_{\frac{a\zeta -b}{d}}(g_{k,b/d}(z))\} .
\eqno(2.14)
$$
We have $u(g_{k,b/d}(z))=u(z)$ and
$$
R^s_{\frac{a\zeta -b}{d}}(g_{k,b/d}(z)) =
(d_{k-\varepsilon_k}g_{k,b/d}^{-1}\left(\frac{a\zeta-b}{d}\right)+b_{k-\varepsilon_k})^{2s}
R^s_{g_{k,b/d}^{-1}\left(\frac{a\zeta-b}{d}\right)}(z)
\eqno(2.15)
$$
This follows from the formula (2.6) on p. 211 of [LZ] and
(2.13). To shorten notation, denote
$$
j_k(b/d,\zeta )^{2s}:=
(d_{k-\varepsilon_k}g_{k,b/d}^{-1}\left(\frac{a\zeta-b}{d}\right)+b_{k-\varepsilon_k})^{2s} .
\eqno(2.16)
$$
Then we get
$$
\lambda_m\psi (\zeta ) =\sum_{ad=m \atop 0< b\le d}\left(\frac{a}{d}\right)^s
\sum_{k=0}^{n(b/d)} (-1)^k j_k(b/d,\zeta )^{2s} \int_{-\infty}^0
 \{u(z), 
R^s_{g_{k,b/d}^{-1}\left(\frac{a\zeta-b}{d}\right)}(z) \} =
$$  
$$
\sum_{ad=m \atop 0< b\le d}\left(\frac{a}{d}\right)^s
\sum_{k=0}^{n(b/d)} (-1)^k j_k(b/d,\zeta )^{2s}
\psi \left(g_{k,b/d}^{-1}\left(\frac{a\zeta-b}{d}\right)\right) .
\eqno (2.17)
$$ 
In order to deduce from (2.17) a nice ``explicit'' formula for $\lambda_m$,
as it was done in [Ma2] for the coefficients of the classical cusp forms,
one could use an appropriate linear functional on functions of $\zeta$.
In the classical case, it was the highest coefficient (or the constant term) of the 
period polynomal.

\smallskip 
 
In the Maass case, one could try to use asymptotic behavious at 0 or $\infty$.
Other forms of Hecke operators, as (2.18) below, might be useful.
 
\medskip

{\bf 3.1. Hecke operators and transfer operator.} T.~M\"uhlenbruch,
using the Choie--Zagier method ([ChZ]),
shows in [M\"uh] that the Hecke operators acting on period
functions for the full modular group can be written in the nice form
$$
T^+_m=\sum_{a>c\ge 0\atop {d>b\ge 0 \atop ad-bc=m}} 
\left( \matrix
a &b\\
c &  d
\endmatrix \right)
\eqno (2.18)
$$
Of course, they act on $\psi (z)$  via $|_s$ in our notation
(M\"hlenbruch denotes this slash operator $|_{2s}$.)

\smallskip

In particular, for $m=1$ we have $T^+_m=I.$

\smallskip

However, if we slightly change the summation domain
replacing $a>c\ge 0$ by $a\ge c>0$, then
then the equations for case $m=1$ will admit the following solutions.
From $ad=1+bc\le 1+(d-1)a$ it follows that $a=c=1$
and $d=b+1\ge 1$ so that we will get the operator
$$
T^*_1:=\sum_{b=0}^{\infty} \left( \matrix
1 &b\\
1 &  b+1
\endmatrix \right)
\eqno (2.19)
$$
This completely {\it ad hoc} correction in fact
makes sense, and moreover, $T^*_1$ imitates
the Hecke operator corresponding to the ``improper
prime $p=1$'', with eigenvalue 1 on $\psi$:

\medskip

{\bf 2.5.1. Claim.} {\it If $\psi (z)$ is a period function
for a Maass cusp form of weight $s$ with $\roman{Re}\,s >0,
s\ne \dfrac{1}{2}$, then}
$$
\psi |_sT^*_1(z)=\psi (z).
\eqno(2.20)
$$

\smallskip

{\bf Proof.} 
Assume moreover that 
$$
\psi^{\tau} (z) := \psi |_s \left( \matrix
0 &1\\
1 & 0
\endmatrix \right) (z) =\varepsilon \psi (z),\ \varepsilon = \pm 1,
\eqno(2.21)
$$
so that $\psi$ is even or odd. This is not a restriction
because any $\psi$ is the sum of an even and odd 
period functions.

\smallskip

According to [LZ2], p. 255, the function 
$$
h(z):=\psi (z+1) =\psi |_s \left( \matrix
1 &1\\
0& 1
\endmatrix \right) (z)
\eqno(2.22)
$$
satisfies the equation
$$
\varepsilon  h|_s \left(\sum_{n=1}^{\infty} 
 \left( \matrix
0 &1\\
1 & n
\endmatrix \right) \right) (z)= \,h(z) .
\eqno(2.23)
$$
Substituting first (2.22) into (2.23), and then (2.21) into resulting identity,
we get:
$$
\psi  |_s \left( \matrix
0 &1\\
1 & 0
\endmatrix \right)     |_s \left( \matrix
1 &1\\
0& 1
\endmatrix \right) |_s \left(\sum_{n=1}^{\infty} 
 \left( \matrix
0 &1\\
1 & n
\endmatrix \right) \right) (z) =
\psi |_s
 \left( \matrix
1 &1\\
0& 1
\endmatrix \right) (z)
\eqno(2.24)
$$
The associativity of the slash operator and the identity
$$
 \left( \matrix
0 &1\\
1 & 0
\endmatrix \right)   \left( \matrix
1 &1\\
0& 1
\endmatrix \right)  
 \left( \matrix
0 &1\\
1 & n
\endmatrix \right)
\left( \matrix
1 &-1\\
0& 1
\endmatrix \right) = \left( \matrix
1 &n-1\\
0& n
\endmatrix \right)
$$
establish (2.20).

\medskip

{\bf 2.6. The Brjuno function and derivatives of the classical $L$--functions.}
The Brjuno function which we referred to in sec. 0.6 is defined as
a generalized L\'evy sum
$$
B(\xi ):= \sum_{j=0}^{\infty} |p_j(\xi )-q_j(\xi )\xi |\,\,\roman{log}\,
\frac{p_{j-1}(\xi )-q_{j-1}(\xi )\xi }{q_j(\xi )\xi -p_j(\xi ) }.
\eqno(2.25)
$$

This series diverges on a set of measure 0.
Outside  it converges to a measurable function, continuous at irrational
points, with period 1. (cf. [MarMouYo2]).

\smallskip

The values of derivatives of  Mellin transforms
of  classical forms were studied by D.~Goldfield ([Go1]) and
N.~Diamantis ([Di]). Goldfeld's  idea consisted in
replacing the $\roman{log}\,y$ initially appearing at the Mellin
expression for the first derivative by the logarithm
of the $\eta$--function, or a combination of such, to
enhance the modular properties of the integrand.
The same game can be played with the Brjuno function
in place of $\eta$--function. 

\smallskip

Consider a classical cusp form $u(z)$ for $SL(2,Z)$ of integral weight
$2k=w+2$ as in 0.1. Let $L_u(s)$ be its Mellin transform.

\medskip

{\bf 2.6.1. Proposition.} {\it   We have  
$$
L_u^{\prime} (w/2+2) =C\,\left[-\int_0^1 u(iy) y^{w/2}B(y)dy
+ \int_1^{\infty} u(iy) y^{w/2-1} B(y)dy \right],
\eqno(2.26)
$$
where}
$$
C= \frac{(2\pi )^{(w+4)/2}}{\Gamma ((w+2)/2)} (1+i^{w+2}) .
$$

\smallskip

{\bf Proof.} First of all, $B(\xi )$ 
satisfies the functional equation
$$
B(\xi ) = -\roman{log}\,\xi +\xi B(\xi^{-1}),\ \xi\in (0,1)
\eqno(2.27)
$$
This is shown by an easy calculation. 

\smallskip

Therefore, we have
$$
\int_0^{\infty} u(iy) y^{w/2} \roman{log}\,y\, dy=
$$
$$
\int_0^{1} u(iy) y^{w/2}\left[ -B(y)+yB(y^{-1})\right] dy+
\int_1^{\infty} u(iv) v^{w/2}\left[ v^{-1}B(v)-B(v^{-1})\right] dv.
$$
In the second summand of the second integrand, make
the change of variable $v=y^{-1}$, and combine it with the
first summand of the first integrand.
Similarly, in the second summand of the first integrand,
make the change of variable  $y=v^{-1}$, and combine it
with the first summand of the second integrand.
This will result in
$$
(1+i^{w+2})\left[-\int_0^1 u(iy) y^{w/2}B(y)dy
+ \int_1^{\infty} u(iy) y^{w/2-1} B(y)dy \right].
$$
The remaining factor in $C$ comes from the Mellin transform.

\medskip

{\it Acknowledgement.} I am grateful to Don Zagier who read
a preliminary draft of this paper and suggested a number of corrections
and complements.

\bigskip

\centerline{\bf References}

\medskip

[ChZ] Y.~Choie, D.~Zagier. {\it Rational period functions for $PSL(2,\Z )$.}
Contemp. Math., 143 (1993), 89--108.

\smallskip

[Di] N.~Diamantis. {\it Special values of higher derivatives of $L$--functions.}
Forum Math., 11:2 (1999), 229--252.

\smallskip

[Go1] D.~Goldfeld. {\it Special values of derivatives of $L$--functions.}
In: Number Theory (Halifax, NS, 1994), CMS. Conf. Proceedings, vol. 15,
159--173, AMS, Providence RI, 1995.

\smallskip

[Go2] D.~Goldfeld. {\it Zeta functions formed with modular symbols.}
In:  Automorphic forms, automorphic representations, and arithmetic (Fort Worth, TX, 1996),  
111--121, Proc. Sympos. Pure Math., 66, Part 1, Amer. Math. Soc., Providence, RI, 1999.

\smallskip

[HiMaMo] J.~Hilgert, D.~Mayer, H.~Movasati. {\it Transfer operator
for $\Gamma_0(n)$ and the Hecke operators for the period
functions of $PSL(2, \bold{Z})$}. Math. Proc. Camb. Phil. Soc.,
139 (2005), 81--116.

\smallskip

[LZ1] J.~Lewis, D.~Zagier. {\it Period functions and
the Selberg zeta function for the modular group.} In:
The Mathematical beauty of physics, Adv. Series in
Math. Phys., 24, World. Sci. Publ., River Edge, NJ, 1997, 83--97.

\smallskip

[LZ2]  J.~Lewis, D.~Zagier. {\it Period functions for Maass wave forms. I.}
Ann. of Math., (2) 153:1 (2001), 191--258.

\smallskip

[M] H.~Maass. {\it \"Uber eine neue Art von nichtanalytischen 
automorphen Funktionen und die Bestimmung
Dirichletsche Reihen durch Funktionalgleichungen.}
Math. Ann. 121 (1949), 141--183.

\smallskip

[MaMar1] Yu.~Manin, M.~Marcolli. {\it Continued fractions, modular symbols, and non-commutative geometry.} Selecta math., new ser. 8:3 (2002),
475--521. e--Print math.NT/0102006

\smallskip

[MaMar2]  Yu.~Manin, M.~Marcolli. {\it Modular shadows and the L\'evy--Mellin
$\infty$--adic transform.} e--Print math.NT/0703718

\smallskip

[Ma1] Yu.~Manin. {\it Parabolic points and zeta-functions of modular curves.}
 Russian: Izv. AN SSSR, ser. mat. 36:1 (1972), 19--66. English:
Math. USSR Izvestija, publ. by AMS, vol. 6, No. 1 (1972), 19--64,
and Selected papers, World Scientific, 1996, 202--247.

\smallskip

[Ma2] Yu.~Manin. {\it Periods of parabolic forms and $p$--adic Hecke series.}
Russian: Mat. Sbornik, 92:3 (1973), 378--401. English:
Math. USSR Sbornik, 21:3 (1973), 371--393
and Selected papers, World Scientific, 1996, 268--290.

\smallskip

[Ma3] Yu.~Manin. {\it The values of $p$--adic Hecke series at integer points
of the critical strip.} Math. USSR Sbornik, 22:4 (1974), 631--637.

\smallskip

[Ma4] Yu.~Manin. {\it Iterated integrals of modular forms and noncommutative modular symbols.}
 37 pp. e--Print math.NT/0502576

\smallskip

[MarMouYo1] S.~Marmi, P.~Moussa, J.-C.~Yoccoz. {\it Complex Brjuno functions.} 
Journ. of the AMS, 14:4 (2001), 783--841.

\smallskip

[MarMouYo2] S.~Marmi, P.~Moussa, J.-C.~Yoccoz. {\it Some properties of real and complex Brjuno functions.}
In: Frontiers in number theory, physics, and geometry. I, 601--623,
Springer, Berlin, 2006.

\smallskip

[May1] D.~Mayer. {\it The thermodynamic formalism approach to Selberg's zeta function for $PSL(2,\Z )$.} Bull. AMS, 25 (1991), 55--60.

\smallskip

[May2] D.~Mayer. {\it  Continued fractions and related transformations.} In: Ergodic Theory, Symbolic Dynamics and Hyperbolic Spaces,
ed. T.Bedford et al., Oxford Univ. Press, NY 1991, 175--222.

\smallskip

[Me] L.~Merel. {\it Universal Fourier expansions of modular forms.}
Springer Levture Notes in Math., 1585 (1994), 95--108.

\smallskip

[M\"uh] T.~M\"uhlenbruch. {\it Hecke operators on period functions 
for the full modular group.} Int. Math. Res. Notices, 77 (2004), 4127--4145.

\smallskip

[Sh1] V.~Shokurov. {\it The study of the homology of Kuga
varieties.} Math. USSR Izvestiya, 16:2 (1981), 399--418.

\smallskip

[Sh2] V.~Shokurov. {\it Shimura integrals of cusp forms.}
Math. USSR Izvestiya, 16:3 (1981), 603--646.

\smallskip

[Za] D.~Zagier. {\it Hecke operators and periods of modular forms.}
In: Israel Math. Conf. Proc., vol. 3 (1990), I.~I.~Piatetski--Shapiro
Festschrift, Part II, 321--336.

\enddocument